\documentclass[11pt,reqno,a4paper]{amsart}
\usepackage{amsthm,amsfonts,amsmath,amssymb}
\usepackage[mathscr]{eucal}
\usepackage[pdftex,bookmarks=true]{hyperref}
\usepackage{enumerate}

\newcommand{\OM}{\Omega}

\newcommand{\bcdot}{\boldsymbol{\cdot}}

\newcommand{\C}{\mathbb{C}}

\newcommand{\N}{\mathbb{N}}

\theoremstyle{plain}
\newtheorem{thm}{Theorem}[section]

\newtheorem{result}[thm]{Result}

\theoremstyle{definition}
\newtheorem{defn}[thm]{Definition}

\theoremstyle{remark}

\title[A Quadrature domain not biholomorphic to a balanced domain]{A $1$-point Quadrature domain of order $1$ not biholomorphic to a
balanced domain}

\author[P.~Haridas]{Pranav Haridas}
\address{Department of Mathematics, Indian Institute of Technology Madras, Chennai 600036, India}
\email{pranav.haridas@gmail.com}
\author[J.~Janardhanan]{Jaikrishnan Janardhanan}
\address{Department of Mathematics, Indian Institute of Technology Madras, Chennai 600036, India}
\email{jaikrishnan@iitm.ac.in}
\thanks{Jaikrishnan Janardhanan is supported by a DST-INSPIRE fellowship from the 
Department of Science and Technology, India.}

\begin{document}

\begin{abstract} 
  It is known that if $f: D_1 \to D_2$ is a polynomial biholomorphism
  with polynomial inverse and constant Jacobian then $D_1$ is a $1$-point
  Quadrature domain (the Bergman span contains all holomorphic polynomials) of
  order $1$ whenever $D_2$ is a balanced domain. Bell conjectured that all 
  $1$-point Quadrature domains arise in this manner. In this note, we construct
  a $1$-point Quadrature domain of order $1$ that is not biholomorphic
  to any balanced domain.
\end{abstract}
\maketitle

\section{Introduction}

  In \cite{bell2009}, Bell conjectures that the any $1$-point
  Quadrature domain of order $1$ in 
  $\C^n, n \geq 2$ is biholomorphic to a balanced domain. 
  A positive result to the above conjecture would have provided
  a characterization analogous to the one provided for planar
  quadrature domains in \cite{MR0447589}. Moreover, as balanced
  domains are very special, the study of holomorphic mappings on $1$-point
  Quadrature domains of order $1$ would have been greatly simplified. 
  Alan Legg \cite[Proposition~4.4]{Legg2017} has
  constructed a $1$-point Quadrature domain 
  that is not biholomorphic to a balanced domain.
  However, the quadrature identity of the
  domain constructed therein was not of order
  $1$.
  
  \medskip

  A large class of 
  domains of interest in several
  complex variables -- including balanced domains,
  quasi-balanced domains and domains invariant under
  the action of a subgroup of $U(n)$ -- are $1$-point
  Quadrature  domains of order $1$. Hence, a complete
  answer to Bell's conjecture on whether every $1$-point
  Quadrature domain of order $1$ is biholomorphic
  to a balanced domain is nevertheless of
  considerable interest.  

  \medskip
  
    In this note, we construct an explicit example of a 
  quasi-balanced domain that is not
  biholomorphic to any balanced 
  domain thereby providing a counterexample to Bell's
  conjecture. We begin with relevant definitions.
  
  \medskip

  \begin{defn}
    Let $D \subset \C^n$ be a bounded domain. The \emph{Bergman span associated
    to $D$} is the $\C$-linear span of the collection
    \[
      \{K_D^\alpha(z,a) : \alpha \in \N^n, a \in D \},
    \]
    where
    \[
      K_D^\alpha(z,a) := \frac{\partial^{|\alpha|}}{\partial 
      \overline{w}^\alpha} K_D(z,w)\bigg|_{w = a}.
    \]
    Here $K_D$ is the Bergman kernel of $D$.
  \end{defn}
  
\begin{defn}
  Let $D \subset \C^n$ be a bounded domain. 
  We say that $D$ is a \emph{quadrature domain} if the constant function 
  $1$ is in the Bergman span of $D$ and a 
  Quadrature domain (with a capital `Q') if all holomorphic
  polynomials are in the the Bergman span of $D$.
  We say that a quadrature (resp. Quadrature) domain $D$ is a 
  $1$-point quadrature (resp. Quadrature) domain of order $1$
  if $K_D(z,0) \equiv c$ for some constant $c\in \C$. 
\end{defn}

\medskip

Observe that the class of Quadrature domains (with a capital `Q')
form a subclass of the class of quadrature domains. Also, if $D$ is a
quadrature domain and $f \in H^2(D)$, the Hilbert space of 
square-integrable holomorphic function on $D$, then we
can find finitely many distinct points $q_1,\dots,q_p \in D$, non-negative
integers
$n_1, \dots, n_q$ and complex numbers $c_{j\alpha}$ (here $\alpha$ is a
multi-index) such that
\[
  \int_D f(z) = \sum_{j = 1}^p \sum_{|\alpha| \leq n_j} c_{j\alpha} 
  \frac{\partial^{|\alpha|}}{\partial z^\alpha} f(z)\bigg|_{z = q_j}.
\]
The above identity  holds independent of the choice of $f$ and is known as a
\emph{quadrature identity for $D$}. The identity is an immediate
consequence of the observation
\[
  \frac{\partial^{|\alpha|}}{\partial z^\alpha} f(z)\bigg|_{z = q_j} =
  \big\langle f, K_D^\alpha(\bcdot,q_j) \big\rangle_{L^2(D)}.
\]
 For a
one point
quadrature (resp. Quadrature) domain of order $1$, it is clear that we have a
quadrature identity of the form
\[
  \int_D f = c \cdot f(a),
\]
where $c \in \C$ and $a \in D$. The recent articles \cite{haridas2015quadrature} and 
\cite{Legg2017} develop the theory of quadrature 
domains in higher dimensions.

\begin{defn}
  Let $p_1,\dots,p_n$ be positive integers with $\gcd(p_1,\dots,p_n) = 1$.
  We say that a
domain $D
\subset \C^n$ is \emph{$(p_1,\dots,p_n)$-balanced (quasi-balanced)} if 
\[
  (\lambda^{p_1}z_1,\lambda^{p_2}z_2,\dots,\lambda^{p_n}z_n) \in D \ \ \forall
  |\lambda|
  \leq 1 \ \ \forall (z_1,z_2,\dots,z_n) \in D.
\]
If $p_1 = p_2 = \dots = p_n = 1$ above, then we say $D$ is a balanced domain 
(also known as a complete circular domain in literature).
\end{defn}

It is well-known that quasi-balanced domains are $1$-point Quadrature domains 
of order $1$ that admit a quadrature identity of the form
\[
  \int_D f dV = \textsf{vol}(D) \cdot f(0).
\]
See \cite[Lemma~2.1]{ning2017} for a proof.
  
\medskip

\section{The Example}

\noindent We require the following two results. 

\begin{result}[Kaup, Folgreng 1 in \cite{kaup1970}]
  Let $D_1, D_2 \subset \C^n$ be two bounded quasi-balanced domains that are
  biholomorphic. Then we can find a biholomorphism $f: D_1 \to D_2$ that
  fixes $0$.
\end{result}

\begin{result}[Ning--Zhou, Theorem~2.2 in \cite{ning2017}]
  Let $D_1, D_2 \subset \C^n$ be two bounded quasi-circular domains that
  contain $0$. Let the corresponding weights be $m = (m_1,\dots,m_n)$ and $m' = 
  (m_1',\dots,m_n')$. Let $f: D_1 \to D_2$ be a biholomorphism that fixes
  $0$. Then writing $f = (f_1,\dots,f_n)$, we have
  \begin{enumerate}
     \item $f$ is a polynomial mapping.
     \item $\textsf{deg}(f_i) \leq \max\{|\delta| : \delta \in \N^n, \delta
     \bcdot m = \gamma \bcdot m, |\gamma| = |\beta|, \beta \bcdot m' = m_i'\}$.
  \end{enumerate} 
  
\end{result}

\noindent \textbf{Construction of the Example.}

Let $D$ be a $(2,3)$-balanced domain in $\C^2$
that is \emph{not} circular. For instance, we might take
\[
  D := \{(z,w) \in \C^2 : |z|^2 + |w|^2 + |z^3 + w^2|^2 < 3 \}.
\]
To see that $D$ as defined above is not circular, observe that $(-1,1) \in D$
but $(-i,i) \not \in D$. 

Suppose $D$ were biholomorphic to a balanced
domain $\OM$. Then
by Kaup's result we can find a biholomorphism $f:\OM \to D$ that fixes $0$. The
result of Ning--Zhou now implies that $f$ has to be a linear. But linear
mappings take balanced domains to balanced domains and $D$ is not balanced by
construction.

\bibliographystyle{amsalpha} \bibliography{Bell_conjecture_counterexample}
\end{document}